\documentclass{article}
\usepackage{amsfonts}
\usepackage{amssymb}

\def\bbbr{{\rm I\!R}}
\def\qed{\rightline{$\square$}}

\newtheorem{definition}{Definition}
\newtheorem{corollary}{Corollary}
\newtheorem{lemma}{Lemma}
\newtheorem{theorem}{Theorem}
\newtheorem{problem}{Problem}
\newtheorem{assumption}{Assumption}
\newcommand{\proof}{\vspace{1mm}\noindent{\it Proof.}\quad}
\newtheorem{example}{Example}

\usepackage{graphicx}

\title{\LARGE \bf
Consensus for switched networks with unknown but bounded
disturbances}

\author{D. Bauso\thanks{D. Bauso is with DINFO, Universit\`a di Palermo, Italy
        {\tt\small dario.bauso@unipa.it}}, L. Giarr\'e\thanks{L. Giarr\'e is with DIAS, Universit\`a di Palermo, Italy
        {\tt\small giarre@unipa.it}}, R. Pesenti\thanks{R. Pesenti is with DMA, Universit\`a ``Ca' Foscari'' di Venezia, Italy
        {\tt\small pesenti@unive.it}} 
\thanks{Research supported by PRIN ``Robustness optimization techniques for high performance
control systems'', and PRIN ``Analysis, optimization, and
coordination of logistic and production systems''.}
}

\begin{document}

\maketitle
\thispagestyle{empty}
\pagestyle{empty}


\small{
\begin{abstract}
We consider stationary consensus protocols for networks of dynamic
agents with switching topologies. The measure of the neighbors'
state is affected by Unknown But Bounded disturbances. Here the
main contribution is the formulation and solution of what we call
the $\epsilon$-consensus problem, where the states are required to
converge in a tube of ray $\epsilon$ asymptotically or in finite
time.
\end{abstract}


}

\section{Introduction.}

\textit{Consensus protocols} are distributed control policies
based on neighbors' state feedback that allow the coordination of
multi-agent systems. According to the usual meaning of consensus,
the system state must converge to an equilibrium point with all
equal components in finite time or
asymptotically~\cite{BGP06,FM02,JLM03,M03,SM04,RBA05,RB05,TJP05,XB04}.

The novelty of our approach is in the presence of Unknown But
Bounded (UBB) disturbances \cite{BR71} in the neighbors' state
feedback. Actually, despite the literature on consensus is now
becoming extensive, only  few approaches have considered a
disturbance affecting the measurements. In our approach we have
assumed an UBB noise, because it requires the least amount of
a-priori knowledge on the disturbance. Only the knowledge of a
bound on the realization is assumed, and no statistical properties
need  to be satisfied.  Moreover, we recall that starting from
[2],  the UBB framework has been used in many different fields and
applications, such as, mobile robotics, vision, multi-inventory,
data-fusion and UAV's and  in estimation, filtering,
identification and robust control theory.

Because of the presence of UBB disturbances convergence to
equilibria with all equal components is, in general, not possible.
The main contribution is then the introduction and solution of the
$\epsilon$-consensus problem, where the states converge in a tube
of ray $\epsilon$ asymptotically or in finite time. In solving the
$\epsilon$-consensus problem we focus on linear protocols and
present a rule for estimating the average from a compact set of
candidate points, say it lazy rule, such that the optimal estimate
for the $i$th agent is the one which minimizes the distance from
$x_i$.

The system under consideration consists of $n$ dynamic agents that
reach consensus on a group decision value by implementing
\emph{distributed} and \emph{stationary control policies} based on
disturbed neighbors' state feedback. Neighborhood relations are
defined by the existence of communication links between nearby
agents. Here, we assume that the set of communication links are
bidirectional and define a time-varying connected communication
network.

The presentation of the results is organized as follows. We first
solve the $\epsilon$-consensus asymptotically and in finite time
for networks with fixed topology (we look at them as switched
systems with dwell time of length infinite). To be more precise,
for a given protocol, we find a tube of minimum radius that the
agents reach asymptotically. Trivially, any tube of radius
strictly greater than the minimum one can be reached in finite
time. We do this by introducing polyhedra of equilibrium points
and studying their stability. The above result means that, in
general, the value of $\epsilon$ cannot be chosen arbitrarily
small. We point out its relation with the amplitude $\xi$ of the
disturbances. We also consider additional assumptions on the
disturbance realization, beyond its inclusion in $D$, and show
that different type of disturbances lead to different values for
the minimum radius. For certain disturbance realizations the
agents are shown to asymptotically reach $0$-consensus. The last
part of this paper, extends the above results to the case of
switching topology. For a given dwell time, we find a tube that
can be reached in finite time. Higher dwell times imply the
convergence to tubes of lower radius.

The paper is organized as follows. In Section~\ref{sec:CP}, we set
up the new framework of switched networks under UBB disturbances
and formulate the $\epsilon$-consensus problem
(Problem~\ref{prob:eps-consensus}). In Section
\ref{sec:LinearProtEqPoints} we introduce the linear protocol and
the lazy rule. In Section~\ref{sec:NonSwitching} we study networks
with fixed topology. In Section \ref{sec:Switching}, we extend the
obtained results to networks with switching topology. Finally, in
Section~\ref{sec:conclusions}, we draw some conclusions.

\section{Switched networks.}\label{sec:CP}
Consider a system of $n$ dynamic agents $\Gamma = \{1, \ldots,
n\}$ and let $\mathcal{E}$ be a finite set of possible edgesets
connecting the agents in $\Gamma$. We model the interaction
topology among agents through a network (graph) $G_{\sigma(t)} =
(\Gamma,E_{\sigma(t)})$, with time variant edgeset $E_{\sigma(t)}
\in \mathcal{E}$, where $\sigma(t)$ is a \emph{switching function}
$\sigma:\bbbr^+ \rightarrow \mathcal{I}$ and $\mathcal{I}$ is the
index set associated with the elements of $\mathcal{E}$. Also, let
us call \emph{switching time}  a time $t$ such that $\sigma(t^-)
\not = \sigma(t^+)$ and let us call \emph{switching interval} the
time interval between two consecutive switching times. In the rest
of this section, to avoid pathological behaviors arising when the
switching times have a finite accumulation point (see, e.g., the
\emph{zeno behavior} in \cite{L03}) and in accordance with
\cite{JLM03,VL05}, we make the following assumption (see, e.g.,
the notion of \emph{dwell time} in~\cite{JLM03,GC05}).

\begin{assumption}\label{asm1} The switching intervals have a finite minimum
length $\tau>0$.\end{assumption}

Henceforth $\tau$ is referred to as the \emph{dwell time}. We also
assume that the edgesets in $\cal E$ induce undirected connected
not complete graphs on $\Gamma$. For each $k \in \mathcal{I}$, the
network $G_k=(\Gamma,E_k)$ is undirected if $(i,j)\in E_k$ then
$(j,i)\in E_k$. The network $G_k$ is connected if for any agent
$i\in \Gamma$ there exists a path, i.e., a sequence of edges in
$E_k$, $(i,i_1)(i_1,i_2)\ldots(i_r,j)$, that connects it with any
other agent $j\in \Gamma$. Finally, the network $G_k$ is not
complete if each agent $i$ is connected (with one edge) only to a
subset of other vertices $N_{ik}=\{j: (i,j) \in E_k\}$ called
\textit{neighborhood of $i$}.\\
Each edge $(i,j)$ in the edgeset $E_k$ means that there is
communication from $j$ to $i$. As $(j,i)$ is also in the edgeset
$E_k$ the communication is bidirectional, namely, if agent $i$ can
receive information from agent $j$ then also agent $j$ can receive
information from agent $i$. Also, $G_k$ not complete means that
each agent $i$ exchanges information only with its neighbors. Here
and in the following, $\mathbf{1}$ stands for the vector
$(1,1,\ldots,1)^T$.

\subsection{Unknown But Bounded disturbances.}\label{sec:UBB}

Let ${\mathcal{T}}$ be the set of switching times. For all $i \in
\Gamma$, consider the family of \emph{first-order} dynamical
systems controlled by a \emph{distributed} and \emph{stationary}
control policy
\begin{eqnarray}
\dot x_i & = & u_{i\sigma(t)}(x_i,y^{(i)}) \quad \forall t \geq 0, ~t \not\in \mathcal{T} \nonumber \\
x_i(t^+) & = & x_i(t^-) \quad \forall t \in
\mathcal{T}\label{eq:xlud}
\end{eqnarray}
where $y^{(i)}$ is the information vector from the agents in
$N_{i\sigma(t)}$  with generic component $j$ defined as follows,
$$ y^{(i)}_j = \left\{
\begin{array}{ll} y_{ij} & \mathrm{if}~  j\in N_{i\sigma(t)}, \\    0 &
\mathrm{otherwise}.
\end{array} \right.$$
In the above equation, $y_{ij}$ is a disturbed measure of $x_j$ obtained by agent $i$ as
$$y_{ij}=x_j+d_{ij}$$
and $d_{ij}$ is an UBB disturbance, i.e., $-\xi \leq d_{ij} \leq
\xi$ with a-priori known $\xi>0$. Hereafter, we denote by $d =
\{d_{ij}, ~(i,j) \in \Gamma^2\}$ the disturbance vector and by $D$
the hypercube $D = \{d: -\xi \leq d_{ij} \leq \xi,~\forall (i,j)
\in \Gamma^2\}$ of the possible disturbance vectors. We assume
that any disturbance realization $\{d(t) \in D, t \geq 0\}$ is
continuous over time.  Note that both $d$ and $D$ are independent
of the topology of network $G$ (which may change over time) as
they are defined on all the possible pairs of agents in $\Gamma$
and not only on the links between them. The continuity hypothesis
on the disturbance realizations can be weakened and most of our
results keep holding true. However, we hold the continuity
assumption to make the proofs of our results simpler and more
readable.

%

\subsection{Problem formulation.}
Before stating the problem we need to introduce the notions of
equilibrium point for a given disturbance realization $d(t)$, and
of $\epsilon$-consensus.

\begin{definition}
A point $x^*$ is an \emph{equilibrium point}
for a given disturbance realization $\{d(t) \in D, t \geq 0\}$ if
there exists $\bar t \geq 0$ such that $u_{\sigma(t)}(x^*_i,
y^{(i)}) = 0$, for all $i \in \Gamma$, for all $t \geq \bar t$.
\end{definition}

According to the usual meaning of consensus, the system state must
converge to an equilibrium point $x^*\in\{\pi \mathbf 1\}$ in
finite time or asymptotically. Hereafter, when we refer to points
of type $\pi \mathbf{1}$, we always understand that $\pi$ may
assume any value in $\mathbb{R}$ and we denote by $\{\pi
\mathbf{1}\}$ the set $\{x:\exists \pi \in \mathbb{R}~s.t.~\pi
\mathbf{1}\}$.

Because of the presence of UBB disturbances convergence to $\{\pi
\mathbf 1\}$ is, in general, not possible. This motivates the
following definition of $\epsilon$-consensus, describing the cases
where the system state is driven in finite time within a bounded
tube of radius $\epsilon$,
\begin{equation}\label{eq:T} T=\left\{x \in
\mathbb{R}^n:\,\left|x_i- x_j\right|\leq 2\epsilon, \, \forall \,
i,j\in \Gamma\right\}.\end{equation}

\begin{definition}
We say that a protocol $u_{\sigma(t)}(.)$ makes the agents reach
\emph{$\epsilon$-consensus in finite time} if there exists a
finite time $\bar t
>0$ such that the system state $x(t)\in T$ for all $t \geq \bar
t$. Furthermore, we say that a protocol $u_{\sigma(t)}(.)$ makes
the agents reach \emph{$\epsilon$-consensus asymptotically}, if
the system state $x(t)\rightarrow T$ for $t \rightarrow \infty$.
\end{definition}

The above definition for $\epsilon =0$ (say it $0$-consensus)
coincides with the usual definition of (asymptotical) consensus.
However, for a generic $\epsilon>0$, the $\epsilon$-consensus in
finite time does not necessarily implies the convergence of the
state $x$ to an equilibrium $x^* \in T$. In other words, $x$ can
be driven to $T$ and keep on oscillating within it for the rest of
the time.

\begin{problem}\label{prob:eps-consensus} ($\epsilon$-consensus problem)
Given the switched system (\ref{eq:xlud}), determine a
(distributed stationary) protocol $u_{\sigma(t)}(.)$ that makes
the agents reach $\epsilon$-consensus in finite time or
asymptotically for any initial state~$x(0)$. Furthermore, study
the dependence of the tube radius $\epsilon$ on the sets
$\mathcal{E}$ and $D$ and on the dwell time $\tau$.
\end{problem}

In the rest of this paper we focus on linear protocols, and
present a rule for estimating the average from a compact set of
candidate points, say it lazy rule, such that the optimal estimate
for the $i$th agent is the one which minimizes the distance from
$x_i$.

\section{Linear protocols and \emph{lazy}
rule.}\label{sec:LinearProtEqPoints}

A typical consensus problem is the average consensus one, i.e.,
the system state converges to the average of the initial state.
Its success derives from the fact that, in absence of
disturbances, it can be simply solved by linear protocols.

Let the linear protocol be given as
\begin{equation}\label{eq:lp}
u_{i\sigma(t)}(x_i,y^{(i)}) =\sum_{j\in N_{i\sigma(t)}} (\tilde
y_{ij}-x_i),  \quad \mbox{for all $i \in \Gamma$}\end{equation}
where $\tilde  y_{ij}$ is the estimate of state $x_j$ on the part
of agent $i$. For a given disturbed measure $y_{ij}$ the state
$x_j$ and consequently its estimate $\tilde y_{ij}$ must belong to
the interval
\begin{equation}\label{eq:tildeyj}
\tilde  y_{ij}\in [ y_{ij}-\xi, y_{ij}+\xi].
\end{equation}

The crucial point is how to select $\tilde  y_{ij}$ from the
interval  $[y_j-\xi,y_j+\xi]$. The next example shows that there
may not exist equilibria if we choose simply $\tilde
 y_{ij}= y_{ij}$.
\begin{example} A three-agent network with a fixed edgeset $E_k$, $N_{1k}=\{1,2\}$, $N_{2k}=\{1,2,3\}$ and
$N_{3k}=\{2,3\}$. A simple criterion is to let $\tilde
 y_{ij}= y_{ij}=x_i+ d_{ij}$
\begin{eqnarray*}\dot x_1 & = &(x_2+d_{12})-x_1\\\dot x_2 & = &[(x_1+d_{21})-x_2] + [(x_3+d_{23})-x_2]\\
\dot x_3 & = &(x_2+d_{32})-x_3 \end{eqnarray*}

Find equilibria by imposing $\dot x=0$ and obtain
\begin{eqnarray}x_1 & = &x_2+d_{12} \\0 & = &d_{12}+d_{21}+d_{23}+d_{32}\\
x_3 & = &x_2+d_{32} \end{eqnarray}

There exist equilibria only if $d_{12}+d_{21}+d_{23}+d_{32}=0,$
that is,  for generic  values of $d_{12}$, $d_{21}$, $d_{23}$, and
$d_{32}$ we cannot guarantee the convergence of the system.
 \end{example}

Let  $\tilde y^{(i)} = \{ \tilde y_{ij}, ~i \in N_i \}$ be defined
according to the \emph{lazy} rule
\begin{equation}\label{eq:tildeyj-1}
\tilde y^{(i)}= \arg \min_{\tilde y_{ij}\in
[y_{ij}-\xi,y_{ij}+\xi], ~j \in N_i\sigma(t)} |\sum_{j\in
N_{i\sigma(t)}} (\tilde y_{ij}-x_i)|.
\end{equation}

Note that as $u(.,.)$ in protocol (\ref{eq:lp}) depends on
$\sum_{j\in N_{i\sigma(t)}} \tilde y_{ij}$, the existence of
multiple solutions $\tilde y^{(i)}$ for (\ref{eq:tildeyj-1}) is
not an issue. This is clearer if one observes that multiple
solutions induce the same value $\sum_{j\in N_{i\sigma(t)}} \tilde
y_{ij}$ for $u(.,.)$ in protocol (\ref{eq:lp}). Given the lazy
rule (\ref{eq:tildeyj-1}), protocol (\ref{eq:lp}) turns out to
have a feedback structure as, for each $i \in \Gamma$, the
quantity $\sum_{j\in N_{i\sigma(t)}} \tilde y_{ij}$ can be
computed as
\begin{equation}\label{eq:sumtildeyj-1} \sum_{j\in N_{i\sigma(t)}}
\tilde y_{ij}= \left\{
\begin{array}{ll}\sum_{j\in N_{i\sigma(t)}}
y_{ij} + |N_{i\sigma(t)}| \xi & if~x_i > \frac{\sum_{j\in
N_{i\sigma(t)}} y_{ij}}{|N_{i\sigma(t)}|}+\xi \\
|N_{i\sigma(t)}| x_i & if~ \frac{\sum_{j\in N_{i\sigma(t)}}
y_{ij}}{|N_{i\sigma(t)}|}-\xi \leq x_i \leq \frac{\sum_{j\in
N_{i\sigma(t)}} y_{ij}}{|N_{i\sigma(t)}|}+\xi\\
\sum_{j\in N_{i\sigma(t)}} y_{ij} - |N_{i\sigma(t)}| \xi & if~x_i
< \frac{\sum_{j\in N_{i\sigma(t)}} y_{ij}}{|N_{i\sigma(t)}|}-\xi
\end{array}\right. .\end{equation}

Hereafter, when we refer to the linear protocol (\ref{eq:lp}), we
always understand that the agents choose $\tilde y^{(i)}$ as in
(\ref{eq:tildeyj-1}).

\section{Fixed topology.}\label{sec:NonSwitching}
In this section we consider a network with fixed topology, i.e., a
network $G = (\Gamma,E)$, with edgeset $E$ constant over time. As
the edgeset $E$ remains constant, for the easy of notation, we
drop the index $\sigma(t)$ from all the notation used throughout
this section. Also when we refer to system (\ref{eq:xlud}) and to
a protocol (\ref{eq:lp}) we always mean that they are associated
to the network $G$.

\subsection{Equilibrium points.}
For a network with fixed topology, we prove that the equilibrium
points exist and belong to polyhedra depending on the type of
disturbance realization. In particular, we state a first result in
the case of constant disturbance $d$, and extend such a result to
the case where the disturbance $d$ takes on values in specific
subsets of $D$.

\begin{lemma}\label{lem:EtaUd}
Given the system (\ref{eq:xlud}) on $G = (\Gamma,E)$, implement a
distributed and stationary protocol $u(.)$ whose components have
the feedback form~(\ref{eq:lp}). If the disturbance $d$ is
constant over time, then:

\begin{itemize}
\item[(i)]  a point $x$ is an equilibrium point for
$u(.)$ if and only if it belongs to the polyhedron

\begin{eqnarray}\label{eq:EquilibriumPoint}P(d,E) &=&
\left\{x: -\frac{\sum_{j\in N_{i}} d_{ij}}{|N_{i}|}- \xi \leq
 \frac{\sum_{j\in N_{i}} x_j }{|N_{i}|} - x_i \leq
-\frac{\sum_{j\in N_{i}} d_{ij}}{|N_{i}|} + \xi, ~ \forall i \in
\Gamma \right\} ;\end{eqnarray}

\item[(ii)] $P(d,E)$ includes all the points in $\{\pi \mathbf
1\}$; in addition, $\{\pi \mathbf 1\}=\bigcap_{d \in D} P(d,E)$.

\item[(iii)]  $P(d,E)$ has $\mathbf 1$ as
only extreme ray up to multiplication by a non-zero scalar.

\end{itemize}

\end{lemma}

\proof{ \emph{(i)} A point $x$ is an equilibrium point if and only
if $u_{i}(x_i,y^{(i)})  = 0$ for all $i \in \Gamma$. This
condition is equivalent to $\min_{\tilde y_{ij}\in
[y_{ij}-\xi,y_{ij}+\xi], ~j \in N_i} |\sum_{j\in N_{i}} (\tilde
y_{ij}-x_i)| = 0$ that, as $y_{ij} =x_j+d_{ij}$, in turn becomes
\begin{eqnarray}
\label{eq:EquilibriumPointXY1}\sum_{j\in N_{i}} (\tilde y_{ij}-x_i) & = & 0, \quad \forall i\in \Gamma \\
\label{eq:EquilibriumPointXY}x_j+d_{ij}-\xi~ \leq~ \tilde y_{ij} &
\leq & x_j+d_{ij}+\xi, \quad  \forall i\in \Gamma, ~\forall j \in
N_i.
\end{eqnarray}
The polyhedron $P(d,E)$ is the projection of the solutions
$(x_i,\tilde y^{(i)})$, for $i \in \Gamma$, of system
(\ref{eq:EquilibriumPointXY1})-(\ref{eq:EquilibriumPointXY}) in
the space of the $x$ variables.

\emph{(ii)} For any $x \in \{ \pi \mathbf 1\}$ it holds that
$\frac{\sum_{j\in N_{i}} x_j }{|N_{i}|} - x_i = 0$. Also,
$-\frac{\sum_{j\in N_{i}} d_{ij}}{|N_{i}|}- \xi \leq 0 \leq
-\frac{\sum_{j\in N_{i}} d_{ij}}{|N_{i}|} + \xi$ because $-\xi
\leq d_{ij} \leq \xi$ for any $i \in \Gamma$, $j \in N_{i}$. Then,
$\{\pi \mathbf 1\} \subseteq P(d,E)$. To prove that $\{\pi \mathbf
1\}=\bigcap_{d \in D} P(d,E)$ we show that $P(\xi,E) = \{\pi
\mathbf 1\}$. To see this last argument, from
(\ref{eq:EquilibriumPoint}) with $d_{ij}=\xi$ for all $i$ and $j$
we have that $\frac{\sum_{j\in N_{i}} x_j }{|N_{i}|} - x_i \leq 0$
for all $i \in \Gamma$ and for any $x \in P(\xi,E)$. The latter
means that the state $x_i$ of each agent $i$ must not be less than
the average state of its neighbors in $N_{i}$ and this situation
occurs only if all the agents have the same state.

\emph{(iii)} The vector $\mathbf 1$ is an extreme ray as it is
immediate to verify that if $x \in P(d,E)$ then  $ x + \pi \mathbf
1 \in P(d,E)$ for any $\pi \in \mathbb{R}$. To prove that a vector
$\mathbf 1$ is the unique extreme ray, up to multiplication by a
non-zero scalar, consider a vector $v$ not parallel to $\mathbf
1$. We note that $0 \in P(d,E)$ and we prove that for some $\pi
\in \mathbb{R}$ the point $0+\pi v \not \in P(d,E)$. As
$-\frac{\sum_{j\in N_{i}} d_{ij}}{|N_{i}|}- \xi$ and
$-\frac{\sum_{j\in N_{i}} d_{ij}}{|N_{i}|} + \xi$ are fixed
values, we have that $\pi v \in P(d,E)$ for any $\pi \in
\mathbb{R}$ if and only if $v_i - \frac{\sum_{j\in N_{i}} v_j
}{|N_{i}|} = 0$ for all $i\in \Gamma$. The latter conditions
define a linear system with $n-1$ independent conditions (provided
that $G$ is connected) and the solutions are of type $v = \eta
\mathbf 1$ for $\eta \in \mathbb{R}$ contradicting the hypothesis
that $v$ is not parallel to $\mathbf 1$.

 \qed}

In the proof of the previous theorem, we have observed that
$-\frac{\sum_{j\in N_{i}} d_{ij}}{|N_{i}|}- \xi \leq 0 \leq
-\frac{\sum_{j\in N_{i}} d_{ij}}{|N_{i}|} + \xi$ for all $i \in
\Gamma$. When such inequalities hold strictly, $P(d,E)$ is a
full-dimensional polyhedron. Actually, any $x$ of type
$(0,\ldots,0, \delta, 0,\ldots,0)$ belongs to $P(d,E)$ if we
choose $\delta > 0$ sufficiently small. However, not all the
polyhedra $P(d,E)$ are full-dimensional as it is apparent by
reminding that $P(\xi,E)=\{\pi \mathbf 1\}$.

In the following we generalize the results of
Lemma~\ref{lem:EtaUd} to the case in which the disturbance is not
constant over time. In other words, we are concerned with the
study of the equilibrium points for generic disturbance
realizations $\{d(t) \in D, ~ t\geq0\}$. First, we can say that
only the points in $\{\pi \mathbf 1\}$ are equilibrium points for
all the possible disturbance realizations $\{d(t) \in D, ~
t\geq0\}$. To see this, observe that i) they are the only
equilibrium points if $d(t)= \xi$ for all $t$ and ii) condition
(\ref{eq:tildeyj-1}) implies $u(x_i, y^{(i)}) = 0$, for all $i \in
\Gamma$ for any realization $\{d(t) \in D, ~ t\geq0\}$.

We will show in the next lemma that, under certain assumptions,
all equilibrium points belong to $P(Q,E)=\bigcup_{d\in Q} P(d,E)$,
for any subset $Q \subseteq D$. Before introducing the lemma,
consider, without loss of generality, the box $Q=\{d \in D:d^-\leq
d \leq d^+\} \subseteq D$ where $d^-$ and $d^+$ are in $D$ and
$d^- \leq d^+$ componentwise. Then it holds
\begin{eqnarray}\label{eq:EquilibriumPointQ}P(Q,E)  &=&
\left\{x\in \mathbb R^n: - \frac{\sum_{j\in N_{i}} d^+_{ij}
}{|N_{i}|}-\xi \leq \frac{\sum_{j\in N_{i}} x_j }{|N_{i}|} - x_i
\leq - \frac{\sum_{j\in N_{i}} d^-_{ij} }{|N_{i}|}+ \xi, ~ \forall
i \in \Gamma \right\}.
\end{eqnarray} To prove (\ref{eq:EquilibriumPointQ}), denote by $\Xi$ the set on the rhs of
(\ref{eq:EquilibriumPointQ}) and note that it holds either $\Xi
\supseteq \bigcup_{d\in Q} P(d,E)$ and $\Xi \subseteq
\bigcup_{d\in Q} P(d,E)$. Actually, for any $d \in Q$, it holds
$d^- \leq d \leq d^+$ then $P(d,E) \subset \Xi$, hence $\Xi
\supseteq \bigcup_{d\in Q} P(d,E)$. Also, to prove $\Xi \subseteq
\bigcup_{d\in Q} P(d,E)$, consider a generic point $\hat x \in
\Xi$. It belongs to $P(\hat d,E)$, where for any $i \in \Gamma$ we
set
\begin{equation}\label{eq:hatd}\hat d_{ij} = \left\{
\begin{array}{ll} d^-_{ij} & if~ j \in N_{i},
~\frac{\sum_{j\in N_{i}} \hat x_j }{|N_{i}|} -
\hat x_i \geq 0\\
d^+_{ij} & if~ j \in N_{i}, ~\frac{\sum_{j\in N_{i}} \hat x_j
}{|N_{i}|} -
\hat x_i < 0\\
d^+_{ij} & otherwise
\end{array} \right. .\end{equation}
As $\hat d \in Q$ by construction, we have $P(\hat d,E) \subseteq
\bigcup_{d\in Q} P(d,E)$ which implies $\Xi \subseteq
\bigcup_{d\in Q} P(d,E)$. Then, we can conclude
that~(\ref{eq:EquilibriumPointQ}) holds true.

In particular, it holds
\begin{eqnarray}\label{eq:EquilibriumPointD}P(D,E)  &=&
\left\{ -2\xi \leq \frac{\sum_{j\in N_{i}} x_j }{|N_{i}|} - x_i
\leq 2\xi, ~ \forall i \in \Gamma \right\}. \end{eqnarray}

Let us define, for a given realization $d(t)$ and a subset $Q$ of
$D$, the value
\begin{equation}\label{eq:mu}
\mu(Q,t_1,t_2) = \max\{\Delta: t_1 \leq \tilde t\leq \tilde t
+\Delta \leq t_2 ~s.t.~ d(t) \in Q~for~all~\tilde t \leq t \leq
\tilde t+\Delta \}.
\end{equation}
In other words, given a time interval $[t_1,t_2]$, the value
$\mu(Q, t_1,t_2)$ is the length of the longest subinterval where
$d(t)$ remains in $Q$.

\begin{lemma}\label{lem:EtaUdIntersezione}
Given the system (\ref{eq:xlud}) on $G = (\Gamma,E)$, implement a
distributed and stationary protocol $u(.)$ whose components have
the feedback form~(\ref{eq:lp}).  Consider a disturbance
realization  $\{d(t) \in D, ~ t\geq0\}$ and box $Q=\{d \in
D:d^-\leq d \leq d^+\} \subseteq D$. Assume that there exist two
nonnegative finite numbers $M$ and $\delta$, such that $\mu(Q,
t,t+M) > \delta$, for all $t \geq 0$. Then, equilibrium points $x$
exist and belong to $P(Q,E)$.
\end{lemma}

\proof{ We first observe that the points $\{\pi \mathbf{1}\}$ are
equilibrium points for a given disturbance realization $d(t)$ and
also that they belong to $P(Q,E)$. We then prove by contradiction
that $x \not \in P(Q,E)$ cannot be an equilibrium point. If $x
\not \in P(Q,E)$, at least for one of its component, say it $i$,
it holds that either $ \frac{\sum_{j\in N_{i}} x_j }{|N_{i}|} -
x_i < - \frac{\sum_{j\in N_{i}} d^+_{ij} }{|N_{i}|}-\xi$ or
$\frac{\sum_{j\in N_{i}} x_j }{|N_{i}|} - x_i > - \frac{\sum_{j\in
N_{i}} d^-_{ij} }{|N_{i}|}+ \xi$. The previous conditions imply
that the value of $u(x_i, y({i}))$ is either strictly less than
zero or strictly greater than zero for all $d \in Q$. Then, for
any $\bar t \geq 0$, there exists a time interval of length
greater than or equal to $\delta$ such that $u(x,y^{(i)})$ is
always either strictly greater than 0 or less than 0. Hence $x$ is
not an equilibrium point.

\qed}

An immediate consequence of the above lemma is the following
corollary.

\begin{corollary}\label{cor:EtaUdIntersezioneMultipla}
Given the system (\ref{eq:xlud}) on $G = (\Gamma,E)$, implement a
distributed and stationary protocol $u(.)$ whose components have
the feedback form~(\ref{eq:lp}).  Consider a disturbance
realization  $\{d(t) \in D, ~ t\geq 0\}$ and a finite set ${\cal
Q} = \{Q_1, Q_2, \ldots\}$ of boxes of $D$. Assume that there
exist two nonnegative finite numbers $M$ and $\delta$, such that
$\mu(Q_r, t,t+M) > \delta$, for all $Q_r \in {\cal Q}$, for all $t
\geq 0$. Then, equilibrium points $x$ for $u(.)$ exist and belong
to $\bigcap_{Q_r \in {\cal Q}} P(Q_r,E)$.
\end{corollary}

In addition, Corollary \ref{cor:EtaUdIntersezioneMultipla} gives
us a hope that if the disturbance realization   enjoys some
general properties the system can reach an equilibrium point close
to the set $\{\pi \mathbf{1}\}$. As an example, consider a
disturbance realization in Corollary
\ref{cor:EtaUdIntersezioneMultipla} characterized, at least, by
${\cal Q} = \{Q_1,Q_2\}$, with $Q_1=\{d \in D:-\xi\leq d \leq
-\hat d\}$ and $Q_2=\{d \in D:\hat d \leq d \leq \xi \}$, with
$0<\hat d \leq \xi$, we obtain that the only equilibrium points
$x$ are in
\begin{eqnarray}\label{eq:EquilibriumPointQ1Q2}P(Q_1,E) \cap P(Q_2,E)  &=&
\left\{  \frac{\sum_{j\in N_{i}} \hat d_{ij} }{|N_{i}|}-\xi \leq
\frac{\sum_{j\in N_{i}} x_j }{|N_{i}|} - x_i \leq -
\frac{\sum_{j\in N_{i}} \hat d_{ij} }{|N_{i}|}+ \xi, ~ \forall i
\in \Gamma \right\}.
\end{eqnarray}
The above set obviously defines a neighborhood of the set $\{\pi
\mathbf{1}\}$, as $\{\pi \mathbf{1}\} \subseteq P(Q,E)$ for any
possible subset $Q$ of $D$.  Interesting is that the radius of the
neighborhood becomes smaller and smaller as $\hat d \rightarrow
\xi$  and that $P(Q_1,E) \cap P(Q_2,E) = \{\pi \mathbf{1}\}$ if
$\hat d = \xi$. The same results hold, for all the situations in
which we can guarantee the disturbance realizations characterized
by ${\cal Q} = \{Q_1,Q_2\}$, such that $P(Q_1,E) \cap P(Q_2,E)$ is
equal to a neighborhood of $\{\pi \mathbf{1}\}$ with a small
radius.

The following corollary asserts that $P(D,E)$ is the minimal set
including all the possible equilibrium points for a policy $u(.)$
given an unknown but bounded disturbance in $D$.

\begin{corollary}\label{cor:EtaUdUnione}
Given the system (\ref{eq:xlud}) on $G = (\Gamma,E)$, implement a
distributed and stationary protocol $u(.)$ whose components have
the feedback form~(\ref{eq:lp}). If the disturbance is unknown but
bounded in $D$, then
\begin{itemize}
\item[(i)] given any $x \in P(D,E)$ then there exists a
disturbance realization $\{d(t) \in D, ~ t\geq0\}$ that has $x$ as
an equilibrium point; \item[(ii)] given any disturbance
realization $\{d(t) \in D, ~ t\geq0\}$ then all its equilibrium
points belong to $P(D,E)$.
\end{itemize}
\end{corollary}

\proof{ \emph{(i)} Any generic point $\hat x \in P(D,E)$ is
trivially an equilibrium point for the corresponding realization
$d(t) = \hat d$, $t\geq 0$, where $\hat d$ is defined as in
(\ref{eq:hatd}).

\emph{(ii)}  Any disturbance realization $\{d(t) \in D:t \leq 0\}$
has, for any finite $M >0$, $\mu(D,t,t+M)=M$ then all its
equilibrium points for $u(.)$ belong to $P(D,E)$ by Lemma
\ref{lem:EtaUdIntersezione}.

\qed}

The following example shows that the value of the parameter
$\epsilon$ defining tube $T$ in (\ref{eq:T}) cannot be chosen
arbitrarily small.

\begin{example} Consider the network $G=(\Gamma,E)$  with
$\Gamma=\{1,\ldots, n\}$ and $E=\{(i,i+1):i=1,\ldots,n-1\}$. Let
$x_{i+1}(t)=x_i(t)+\xi$ for any arbitrary value of $x_1(t)$. This
point is an equilibrium as long as $d_{ij}(t)=0$ for all $i\in
\Gamma$ and $t\geq 0$. In this situation, the value $\epsilon$
defining $T$ in (\ref{eq:T}) is equal to $\frac{n-1}{2}
\xi$.\end{example}

Corollary \ref{cor:EtaUdUnione} suggests a way to determine a
strict upper bound $\bar \epsilon$ for $\epsilon$. We have
\begin{equation}\label{eq:valueOfEpsilon}
\bar \epsilon =  \max_{i,j \in \Gamma} \max_{x \in P(D,E)} \{x_i -
x_j\},
\end{equation}
whose brute force computation requires the solution of $n(n-1)$
linear programming problems of type $\max_{x \in P(D,E)} \{x_i -
x_j\}$. Then, the computation of $\bar \epsilon$ becomes
polynomial.


\subsection{Stability.}\label{sec:Stability}

In this subsection we prove the asymptotic stability of the
equilibrium points. To this end, we have to introduce a basic
property of the stationary protocol $u(.)$ whose components have
the feedback form~(\ref{eq:lp}). In the following, we denote by
$sign: \mathbb{R}\rightarrow \{-1,0,1\}$ the function that returns
1 if its argument is positive, -1 if its argument is negative, 0
if its argument is null.

\begin{lemma}\label{Lem:SignPerm}
Given the system (\ref{eq:xlud}) on $G = (\Gamma,E)$, implement a
distributed and stationary protocol $u(.)$ whose components have
the feedback form~(\ref{eq:lp}).  Either $\mbox{\rm
sign}(u_{i}(x_i,y^{(i)})) = \mbox{\rm sign}(\sum_{j\in N_{i}} (
x_{j}-x_i))$ or $\mbox{\rm sign}(u_{i}(x_i,y^{(i)})) = 0$, for
each $i \in \Gamma$, for each $t \geq 0$.

\end{lemma}
\proof{For each $i \in \Gamma$, for each $t \geq 0$, given the
protocol  $u_{i}(x_i,y^{(i)}) =\sum_{j\in N_{i}} (\tilde
y_{ij}-x_i)$, consider the solution of the linear problem that
defines the value of $\tilde y^{(i)}$
\begin{equation}\label{eq:tildeyj-2}
z_i =  \min_{\tilde y_{ij}\in [y_{ij}-\xi,y_{ij}+\xi], ~j \in N_i}
|\sum_{j\in N_{i}} (\tilde y_{ij}-x_i)|.
\end{equation}
If $z_i = 0$ the lemma is proved. If $z_i > 0$ two situations can
occur,  the value of $\sum_{j\in N_{i}} (\tilde y_{ij}-x_i)$ is
either strictly positive or strictly negative, for any $\tilde
y_{ij}\in [y_{ij}-\xi,y_{ij}+\xi]$, $j \in N_i$.  We claim that if
$\sum_{j\in N_{i}} (\tilde y_{ij}-x_i)
> 0$ then $u_{i}(x_i,y^{(i)}) > 0$ and $\sum_{j\in
N_{i}} ( x_{j}-x_i) > 0$. If $\sum_{j\in N_{i}} (\tilde
y_{ij}-x_i) > 0$ for any $\tilde y_{ij}\in
[y_{ij}-\xi,y_{ij}+\xi]$, $j \in N_i$ then, by definition,
$u_{i}(x_i,y^{(i)}) > 0$, as the chosen $\tilde y_{ij}$ must
belong to $[y_{ij}-\xi,y_{ij}+\xi]$. In addition, we have
$\sum_{j\in N_{i}} (x_{j} + d_{ij} - \xi -x_i) > 0$, hence
$\sum_{j\in N_{i}} (x_{j} - x_i) > \sum_{j\in N_{i}}(\xi - d_{ij})
\geq 0$, as  $-\xi \leq d_{ij} \leq \xi$.

A symmetric argument holds if $\sum_{j\in N_{i}} (\tilde
y_{ij}-x_i) < 0$, for any $\tilde y_{ij}\in
[y_{ij}-\xi,y_{ij}+\xi]$, $j \in N_i$.

\qed}

\begin{theorem}\label{Th:EquilPoint2}
Given the system (\ref{eq:xlud}) on $G = (\Gamma,E)$, implement a
distributed and stationary protocol $u(.)$ whose components have
the feedback form~(\ref{eq:lp}). Then, the system trajectory
converges to equilibrium points in~$P(D,E)$.
\end{theorem}

\proof{ We prove the convergence to equilibrium points in $P(D,E)$
by introducing a candidate Lyapunov function $V(x)=
\frac{1}{2}\sum_{(i,j) \in E}(x_j-x_i)^2$. Trivially, $V(x) = 0$
if and only if $x \in \{\pi \mathbf{1}\}$; $V(x) > 0$ for all $x
\not \in \{\pi \mathbf{1}\}$. We now prove that $\dot V(x) < 0$
for all $x \not \in P(D,E)$. On this purpose, for $\dot V(x)$ we
can write

\begin{eqnarray}
\dot V(x) &= &  \sum_{(i,j) \in E} (x_j-x_i)(u_{j}-u_{i})=
-\sum_{i \in \Gamma}u_{i} \sum_{j \in N_{i}}(x_j-x_i)= \nonumber \\
 &= & -\sum_{i \in \Gamma} \mbox{sign}(u_{i}) \, \mbox{sign}\left(\sum_{j \in
 N_{i}}(x_j-x_i)\right)
\left|u_{i}\right| \left|\sum_{j \in N_{i}}(x_j-x_i)\right|
\label{eq:StabOfxiE}
\end{eqnarray}

From Lemma \ref{Lem:SignPerm}, if $\sum_{j \in N_{i}}(x_j-x_i) =
0$ then $u_{i}(x,y^{(i)}) = 0$. This in turns implies that $\dot
V(x)$ is null if and only if $u_k(x) = 0$. The latter observation
is sufficient to prove that i) the state trajectory converges to
$P(D,E)$ and that ii) the convergence is to an equilibrium point.
Indeed, for $t \rightarrow \infty$, we have $\dot V \rightarrow
0$. Then, $u \rightarrow 0$ and consequently $\dot x \rightarrow
0$.

\qed}

\begin{theorem}\label{Th:EquilPoint3}
Given the system (\ref{eq:xlud}) on $G = (\Gamma,E)$, implement a
distributed and stationary protocol $u(.)$ whose components have
the feedback form~(\ref{eq:lp}).  Consider a disturbance
realization  $\{d(t) \in D, ~ t\geq0\}$ and a box $Q=\{d \in
D:d^-\leq d \leq d^+\} \subseteq D$. Assume that there exist two
nonnegative finite numbers $M$ and $\delta$, such that $\mu(Q,
t,t+M) > \delta$, for all $t \geq 0$. Then, the system trajectory
converges to equilibrium points in  $P(Q,E)$.
\end{theorem}

\proof{ We prove the convergence to equilibrium points in $P(Q,E)$
following the same argument used in the proof of Theorem
\ref{Th:EquilPoint2}. We only note that now we have more
information on the disturbance. In particular we know that, for
every time interval of length $M$, it spends at least a time
$\delta$ assuming values in $Q$. The explicit dependence of the
disturbance on time makes the Lyapunov function time-varying.

To prove the system stability we make use of the results in
\cite{Ma03}. We define the  function $p:
\mathbb{R}\rightarrow[0,\infty)$
\begin{equation}\label{eq:ptdefinition}
p(t) = \left\{\begin{array}{ll} 1 & ~if~d(t) \in Q \\
0 & otherwise \end{array}\right.
\end{equation}
It is immediate to verify that $p(.)$ satisfies the conditions in
Remark 3 in \cite{Ma03}, in particular, there exists three finite
values $\bar p$, $M$, $\delta >0$ such that $0 \leq p(t) \leq \bar
p$, $\int_{t}^{t+M} p(s)ds \geq \delta$ for all $t\geq 0$. We also
define the function $W: \mathbb{R}^n\rightarrow[0,\infty)$
\begin{equation}\label{eq:Wdefinition}
W(x) = \left\{\begin{array}{ll} 0 & ~if~x \in P(Q,E) \\
\min_{\tilde y^{(i)}\in [x+d-\xi,x+d+\xi], ~d \in Q}
\left\{|\sum_{j\in N_{i}} (\tilde y_{ij}-x_i)|\left|\sum_{j \in
N_{i}}(x_j-x_i)\right|\right\} & ~if~x \not \in P(Q,E)
\end{array}\right.
\end{equation}
Observe that $W(x) = 0$ for $x \in P(Q,E)$, whereas $0< W(x) <
\left|u_{i}\right| \left|\sum_{j \in N_{i}}(x_j-x_i)\right|$ for
all $x \not \in P(Q,E)$ and for all $t \geq 0$. Hence $\dot V(x)
\leq -p(t)W(x) \leq 0$ for all $x$ and all $t\geq0$ and, in
particular, $\dot V(x)= 0$ for all $t \geq 0$, only for  $~x \in
P(Q,E)$. The system trajectory converges to $P(Q,E)$. Finally, we
note that the system trajectory converges to an equilibrium point
as $\dot V(x)$ is null if and only if $u_k(x) = 0$.

\qed}

An immediate consequence of the above theorem is the following
corollary.

\begin{corollary}\label{cor:TheoEtaUdIntersezioneMultipla}
Given the system (\ref{eq:xlud}) on $G = (\Gamma,E)$, implement a
distributed and stationary protocol $u(.)$ whose components have
the feedback form~(\ref{eq:lp}).  Consider a disturbance
realization  $\{d(t) \in D, ~ t\geq 0\}$ and a finite set ${\cal
Q} = \{Q_1, Q_2, \ldots\}$ of boxes of $D$. Assume that there
exist two nonnegative finite numbers $M$ and $\delta$, such that
$\mu(Q_r, t,t+M) > \delta$, for all $Q_r \in {\cal Q}$, for all $t
\geq 0$. Then, the system trajectory converges to $\bigcap_{Q_r
\in {\cal Q}} P(Q_r,E)$.
\end{corollary}

Finally, we can conclude that for a disturbance realization that
in Corollary \ref{cor:EtaUdIntersezioneMultipla} and in Corollary
\ref{cor:TheoEtaUdIntersezioneMultipla} is characterized, at
least, by ${\cal Q} = \{Q_1,Q_2\}$, with $Q_1=\{d \in D:-\xi\leq d
\leq -\hat d\}$ and $Q_2=\{d \in D:-\hat d \leq d \leq \xi \}$,
with $0<\hat d \leq \xi$,  the  system trajectory converges to a
neighborhood of the set $\{\pi \mathbf{1}\}$ with the ray of the
neighborhood that becomes smaller and smaller as $\hat d
\rightarrow \xi$ and $P(Q_1,E) \cap P(Q_2,E) = \{\pi \mathbf{1}\}$
if $\hat d = \xi$.


\subsection{Bounds for $x(t)$.}
In this subsection, we determine bounds for the minimum and
maximum value that the components of $x(t)$ assume over the time
depending on the initial state $x(0)$ and for any disturbance
realization $\{d(t)\in D: t\geq 0\}$. In particular, we prove
that, when we apply the lazy rule (\ref{eq:tildeyj-1}), we always
obtain
\begin{equation}
\alpha(x(0)) \leq \lim_{t \rightarrow \infty} x_i(t) \leq
\beta(x(0)), \quad \mbox{ for all $i \in \Gamma$}.
\end{equation}
As a further result, we also show that the difference between the
maximum and the minimum agent states may not increase over the
time.

Given a network $G$ and an initial state $x(0)$, the main idea is
to replace $G$ by a much simpler network $H$ composed by only two
agents and such that the initial state of $H$ is equal to the two
maximal values of the initial state of network $G$. The result is
that the maximal value assumed by the states of $G$ is always
bounded by the values assumed by the states of $H$.

Let us denote by $H = (\{a,b\},\{(a,b)\})$ a system with only two
connected agents. Let $x^H({t})$ be the state of $H$, namely,
$x^H(t)=\{x^H_a(t),x^H_b(t)\}$. Let the components of $x^H({t})$
be subject to a constant disturbance $d^H(t) = \xi$. Let us also
define $i_1(t) =arg\max_{j\in \Gamma}\{x_j(t)\}$ and  $i_2(t) =
arg\max_{j\in \Gamma \setminus \{i_1(t)\}}\{x_j(t)\}$, for all $t
\geq 0$. Actually, $i_1(t)$ and $i_2(t)$ are the two agents with
the first two maximal states. Obviously, $i_1(t)$ and $i_2(t)$
depend on time $t$. Analogously, define $i_n(t) =arg\min_{j\in
\Gamma}\{x_j(t)\}$ and by $i_{n-1}(t) = arg\min_{j\in \Gamma
\setminus \{i_n(t)\}}\{x_j(t)\}$, for all $t \geq 0$.

\begin{lemma}\label{lem:initial state2}
Given the system (\ref{eq:xlud}) on $G = (\Gamma,E)$ with initial
state $x({0})$, implement a distributed and stationary protocol
$u(.)$ whose components have the feedback form~(\ref{eq:lp}).
Consider the system $H = (\{a,b\},\{(a,b)\})$ with an initial
state $x^H_a(0) = x_{i_1(0)}(0)$ and $x^H_b(0) = x_{i_2(0)}(0)$.
Then, for all $t \geq 0$,  $x^H_a(t) \geq x_{i}(t)$, for all $i
\in \Gamma$ and $x^H_b(t) \geq x_{i}(t)$, for all $i \in \Gamma
\setminus \{i_1(t)\}$.
\end{lemma}

\proof {Observe that $x_i(t)$ is a differentiable variable for all
$i \in \Gamma$. The same property holds for $x^H_a(t)$ and
$x^H_b(t)$. In addition it holds that $x^H_a(t) \geq x^H_b(t)$ for
any $t \geq 0$.

At time $t=0$ the thesis holds by definition of values $x^H_a(0)$
and $x^H_b(0)$. By contradiction, assume that at some time instant
$\bar t> 0$ the thesis is false, i.e., there exist some $i,j \in
\Gamma$ such that either $x^H_a(\bar t) < x_{i}(\bar t)$ or
$x^H_a(\bar t) \ge x_{i}(\bar t)$ but $x^H_b(\bar t) < x_{j}(\bar
t)$. By continuity, there must also exists $0\leq  t < \bar t$
where one of the following conditions holds
\begin{itemize}
\item[i)] $x^H_a(t) = x_{i}(t)$, $x^H_b(t) = x_{j}(t)$,
          $x^H_b(t) \geq x_{k}(t)$, for all $k \in \Gamma \setminus\{i,j\}$,
          and either $x^H_a(t+dt) < x_{i}(t+dt)$ or $x^H_b(t+dt) < x_{j}(t+dt)$;
\item[ii)] $x^H_a(t) = x_{i}(t)$,
          $x^H_b(t) > x_{k}(t)$, for all $k \in \Gamma \setminus\{i\}$, and $x^H_a(t+dt) < x_{i}(t+dt)$;
\item[iii)] $x^H_a(t) > x_{i}(t)$, $x^H_b(t) = x_{j}(t)$,
          $x^H_b(t) \geq x_{k}(t)$, for all $k \in \Gamma \setminus\{i,j\}$, and $x^H_b(t+dt) < x_{j}(t+dt)$;
\end{itemize}

Consider case i). It holds $\dot x^H_a(t) =
\tilde{y}^H_{ab}(t)-x^H_a(t) \leq 0$ and $\dot x_i(t) = \sum_{r\in
N_{i}} (\tilde y_{ir}(t)-x_i(t)) \leq 0$ and, in particular,
$\tilde y_{ir}(t)-x_i(t) \leq 0$, for all $r\in N_{i}$. As by
hypothesis $\tilde{y}^H_{ab}(t)-x^H_a(t) \geq \tilde
y_{ir}(t)-x_i(t)$ for any $r \in \Gamma$, we have $\dot
x^H_a(t)\geq \dot x_i(t)$, hence the inequality $x^H_a(t+dt) <
x_{i}(t+dt)$ is false. It also holds $\dot x^H_b(t) =
\tilde{y}^H_{ba}(t)-x^H_b(t) \geq 0$ and $\dot x_j(t) = \sum_{r\in
N_{j\sigma(t)}} (\tilde y_{jr}(t)-x_j(t)) \leq \sum_{r\in
N_{j\sigma(t)}\setminus \{i\}} (\tilde y_{jr}(t)-x_j(t)) + (\tilde
y_{ji}(t)-x_j(t)) $. As by hypothesis
$\tilde{y}^H_{ba}(t)-x^H_b(t) \geq \tilde y_{ji}(t)-x_j(t)$ and
$\sum_{r\in N_{j\sigma(t)}\setminus \{i\}} (\tilde
y_{jr}(t)-x_j(t)) \leq 0$, we have $\dot x^H_a(t)\geq \dot
x_i(t)$, hence the inequality $x^H_b(t+dt) < x_{j}(t+dt)$ is
false. Hence Conditions i) cannot hold.

We can use the first part of the above argument to prove that
Conditions ii) cannot hold, and use the second part to prove that
Conditions iii) cannot hold.

\qed}

The following corollary holds.
\begin{corollary}\label{cor:beta}
Given the system (\ref{eq:xlud}) on $G = (\Gamma,E)$ with initial
state $x({0})$, implement a distributed and stationary protocol
$u(.)$ whose components have the feedback form~(\ref{eq:lp}).
Then:

\begin{itemize}
\item[(i)]

The values assumed by the state trajectory $x(t)$ for $t
\rightarrow \infty$ satisfy the following inequalities for any
disturbance realization $d(t)$

\begin{equation}\label{eq:bounds}
x_{i_n}(0)\leq \alpha(x(0)) \leq \lim_{t \rightarrow \infty} x_i(t) \leq
\beta(x(0))\leq x_{i_1}(0), \quad \mbox{ for all $i \in \Gamma$},
\end{equation}

where the bounds $\alpha(x(0))$ and $\beta(x(0))$ depend on the
initial state $x(0)$ as follows
\begin{eqnarray}
\alpha(x(0)) &  = & \max\left\{x_{i_n}(0), \frac{x_{i_n}(0)+
x_{i_{n-1}}(0)}{2} - \xi -
\frac{\xi}{2}\ln\frac{-x_{i_n}(0)+x_{i_{n-1}}(0)+\xi}{\xi}\right\}\\
\beta(x(0)) & = &\min\left\{x^H_a(0), \frac{x^H_a(0)+ x^H_b(0)}{2}
+ \xi + \frac{\xi}{2}\ln\frac{x^H_a(0)-x^H_b(0)-\xi}{\xi}\right\}.
\end{eqnarray}

\item[(ii)] The value $x_{i_1}(t)$ is non increasing on $t$ and
the value $x_{i_n}(t)$ is non decreasing on $t$.
\end{itemize}

\end{corollary}

\proof{ To prove (i), consider the network $H$ as defined in the
proof of Lemma \ref{lem:initial state2} and let us study the
evolution of $x^H(t)$ over the time. If $x^H_a(0)-x^H_b(0) \leq
2\xi$, the system state evolves according to $\dot x^H_a(t) = 0$
and $\dot x^H_b(t) = x^H_a(0) - x^H_b(t)$.
 If $x^H_a(0)-x^H_b(0)>2\xi$,
the system state evolves according to $\dot x^H_a(t) =  x^H_b(t) +
2\xi - x^H_a(t)$ and $\dot x^H_b(t) = x^H_a(t) - x^H_b(t)$, as
long as $0 \leq t \leq \hat t$, where $\hat t$ is such that $
x^H_a(\hat t) - x^H_b(\hat t) = 2\xi $, i.e., $\hat t =
\frac{1}{2}\ln\frac{x^H_a(0)-x^H_b(0)-\xi}{\xi}$. For $t \geq \hat
t$, the system state evolves according to $\dot x^H_a(t) = 0$ and
$\dot x^H_b(t) = x^H_a(\hat t) - x^H_b(t)$. Hence, for $t \geq 0$,
we have

\begin{eqnarray}\label{eq:xat}
x^H_a(t) &=& \left\{\begin{array}{ll}
        x^H_a(0) &~if~x^H_a(0)-x^H_b(0) \leq 2\xi\\
        \frac{\xi(1+2t-{e^{-2t}})+ x^H_a(0)(1+{e^{-2t}})+ x^H_b(0)(1-{e^{-2t}})}{2}
                 &~if~x^H_a(0)-x^H_b(0) > 2\xi~and~ t\leq \hat t\\
        \frac{x^H_a(0)+ x^H_b(0)}{2} + \xi + \frac{\xi}{2}\ln\frac{x^H_a(0)-x^H_b(0)-\xi}{\xi}
                 &~if~x^H_a(0)-x^H_b(0) > 2\xi~and~ t> \hat t
\end{array}\right.\\
x^H_b(t) &=& \left\{\begin{array}{ll}
        x^H_a(0) - x^H_b(0) e^{-t} &~if~x^H_a(0)-x^H_b(0) \leq 2\xi\\
        \frac{\xi(-1+2t+{e^{-2t}})+ x^H_a(0)( 1-{e^{-2t}})+ x^H_b(0)(1+{e^{-2t}})}{2}
                 &~if~x^H_a(0)-x^H_b(0) > 2\xi~and~ t\leq \hat t\\
        x^H_a(\hat t) - x^H_b(\hat t) e^{-t}
                 &~if~x^H_a(0)-x^H_b(0) > 2\xi~and~ t> \hat t
\end{array}\right. .\label{eq:xbt}
\end{eqnarray}

For $t\rightarrow \infty$, $x^H_a(t)$ and $x^H_b(t)$ converge to
$\beta(x(0)) = \min\left\{x^H_a(0), \frac{x^H_a(0)+ x^H_b(0)}{2} +
\xi + \frac{\xi}{2}\ln\frac{x^H_a(0)-x^H_b(0)-\xi}{\xi}\right\}$.
Hence, $\lim_{t\rightarrow \infty} x_i(t) \leq \beta(x(0))$, for
all $i \in \Gamma$. With an analogous argument, we can prove
$\lim_{t\rightarrow \infty} x_i(t) \geq \alpha(x(0)) =
\max\left\{x_{i_n}(0), \frac{x_{i_n}(0)+ x_{i_{n-1}}(0)}{2} - \xi
-
\frac{\xi}{2}\ln\frac{-x_{i_n}(0)+x_{i_{n-1}}(0)+\xi}{\xi}\right\}$,
for all $i \in \Gamma$.

To prove (ii), observe that we have $$x_{i_1}(0)= x_a^H(0)\geq
x_a^H(t) \geq x_{i_1}(t),$$ where the first equality and the last
inequality hold by definition, whereas the inequality
$x_a^H(0)\geq x_a^H(t)$ derives straightforwardly from the fact
that $\dot x_a^H(t)\leq 0$ for all $t\geq 0$. \qed}

Corollary \ref{cor:beta} (ii) proves that the system trajectory
$x(t)$ is bounded as $t \rightarrow \infty$ and also that the
difference between the maximum and the minimum agent states may
not increase over the time. More formally, denote by
$\mathcal{V}(x(t)) = x_{i_1}(t) - x_{i_n}(t)$ then
\begin{equation}\label{eq:dec}\mathcal{V}(x(t)) \geq \mathcal{V}(x(t+\Delta t))\quad \mbox{for any
$t \geq 0$ and $\Delta t > 0$.}\end{equation} We use this last
implication to introduce some additional results that will turn
useful when dealing with switching topology systems.

Denote by  $\mathcal{V}_\infty = \lim_{t \rightarrow
\infty}(x_{i_1}(t) - x_{i_n}(t))$ the final value of
$\mathcal{V}(x(t))$. Observe that for some network $G(\Gamma,E)$
and initial state $x(0)$, there may exist some disturbance
realizations $\{d(t)\in D: t\geq 0\}$ such that even if
$\mathcal{V}(x(0)) > \mathcal{V}_\infty$, the value
$\mathcal{V}(x(t))$ may be constant over some finite time interval
before reaching its final value $\mathcal{V}_\infty$. More
specifically, there may exist $t$ and $\Delta t$ such that
$\mathcal{V}(x(t)) = \mathcal{V}(x(t+\Delta t)) >
\mathcal{V}_\infty$.
\begin{example}
Consider a networks of six agents with chain topology depicted in
Fig.~\ref{fig:Chain6agents}. The initial state is
$x(0)=[100,100,100,0,0,0]^T$ and disturbances are
$d_{12}=d_{21}=d_{23}=d_{32}=d_{34}=1$ and
$d_{43}=d_{45}=d_{54}=d_{56}=d_{65}=-1$. Figure \ref{fig:Zoom} (a)
shows the time plot of the evolution of the state $x(t)$ for
$0\leq t \leq 20$ and as it can be seen, trajectories converge to
the equilibrium $x^*=[63,61,55,45,39,37]^T$ with $\epsilon = 26$
(note that the initial deviation between maximum and minimum value
of the state is $100$). Figure \ref{fig:Zoom} (b) displays a zoom
of the trajectories for $0\leq t\leq 3$ pointing out that
$\mathcal{V}(x(t)) = \mathcal{V}(x(0))$ for $0\leq t \leq 0.5$.
\begin{figure}[H]
    \centering%
    \includegraphics[scale=1]{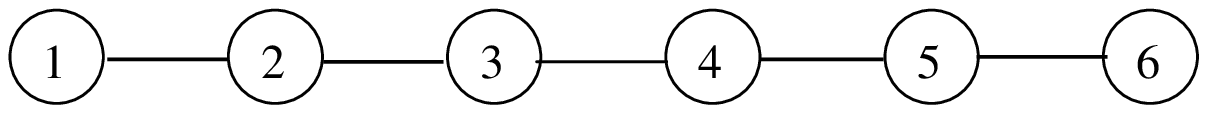}
    \caption{Chain of six agents.}
    \label{fig:Chain6agents}
\end{figure}

\begin{figure}[H]
    \centering%
    \includegraphics[scale=1]{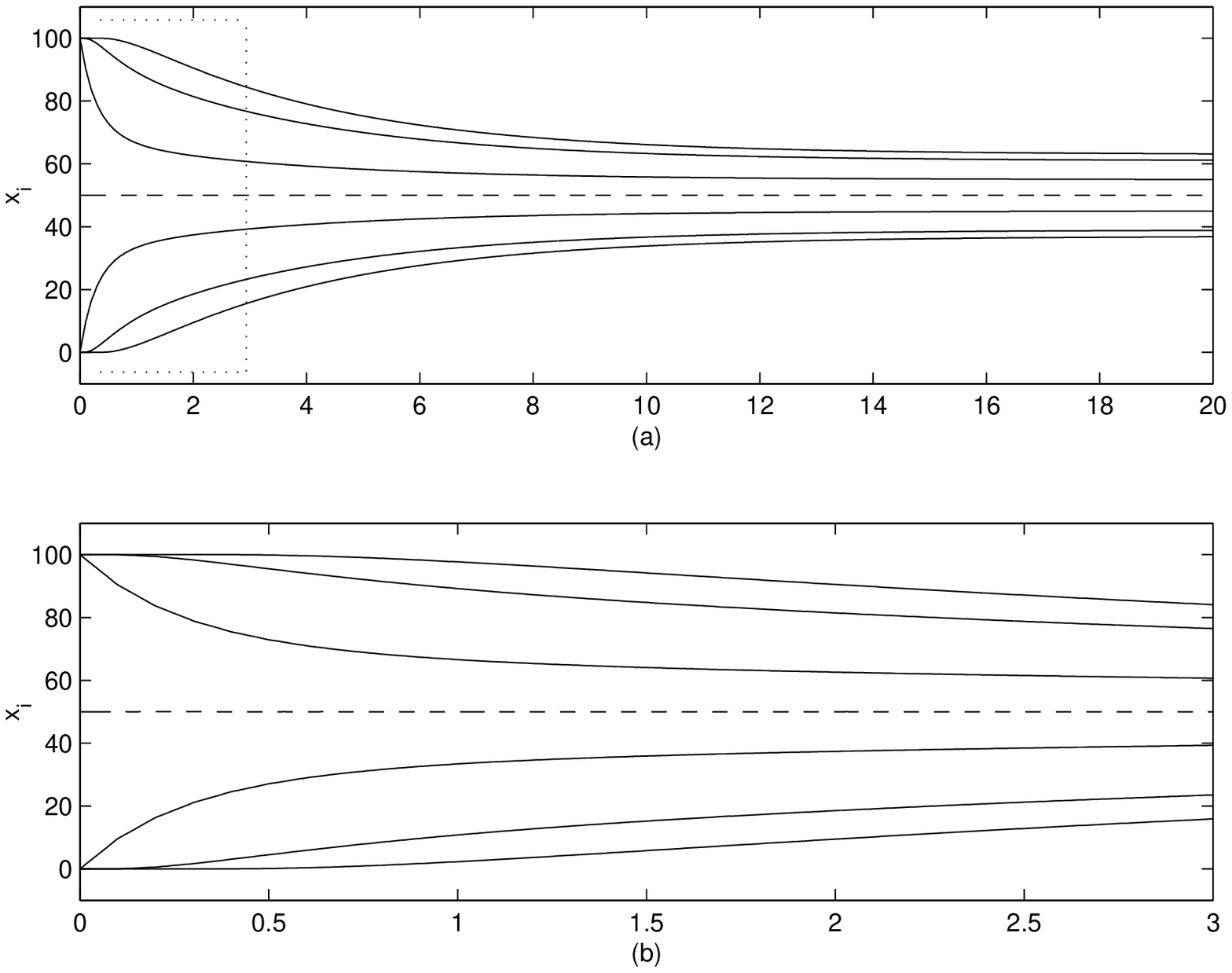}
    \caption{(a) Time plot of the state $x$ for an initial
    state $x(0)=[100,100,100,0,0,0]^T$. Trajectories converge to the equilibrium
    $x^*=[63,61,55,45,39,37]^T$ with $\epsilon = 26$; (b) zoom of the first time instants (dotted rectangle in (a)),
    which highlights $\mathcal{V}(x(t)) = \mathcal{V}(x(0))$ for $0\leq t \leq 0.5$.}
    \label{fig:Zoom}
\end{figure}
\end{example}

In the following we prove that we can always sample the state
trajectory in such a way that the sequence of values for ${\cal
V}(.)$ is strictly decreasing on time until the state is
``almost'' in $P(D,E)$.

 To this end, given the system~(\ref{eq:xlud}) on $G =
(\Gamma,E)$ and disturbance realizations $\{d(t)\in D: t\geq 0\}$,
with a little abuse of notation, we denote by $\mathcal{V}(x) =
\max_{i \in \Gamma}\{x_i\} - \min_{i \in \Gamma}\{x_i\}$ the
maximum difference between two components of vector $x$, for any
$x \in \mathbb{R}^n$. Also, we denote by $P(D,E)+ \nu = \{x \in
\mathbb{R}^n: \exists y \in P(D,E)~s.t.~\|x-y\|_\infty \leq \nu
\}$, with $\nu >0$, the set of points whose distance from set
$P(D,E)$ is not greater than $\nu$, according the
$\mathcal{L}_\infty$ norm; and by $V(x) =\frac{1}{2} \sum_{(i,j)
\in E}(x_i -x_j)^2$ the Lyapunov function considered in the proof
of Theorem \ref{Th:EquilPoint2}.

The following lemma holds
\begin{lemma}\label{lem:VandL}
Given the system (\ref{eq:xlud}) on $G = (\Gamma,E)$. Let $\hat x$
and $\bar x$ in $\mathbb{R}^n$, if $V(\bar x) \leq
\frac{4\gamma^2V(\hat x)}{n^2(n-1)}$ then $\mathcal{V}(\bar x)\leq
{\gamma}\mathcal{V}(\hat x)$, for any $0 < \gamma < 1$.
\end{lemma}
\proof{ First we determine the bounds for the values of $V(x)$ for
$x \in \mathbb{R}^n$ such that $\mathcal{V}(x) = \mathcal{V}=
const$. Denote by $\hat E = \{(i,j): i<j, ~i,j \in \Gamma\}$ the
edgeset of the complete network induced by vertices in $\Gamma$.
Observe that if $x \in \mathbb{R}^n$ and $\mathcal{V}(x) =
\mathcal{V}(x^1)$ then
$$V(x) \leq \max_{x \in \mathbb{R}^n} \left\{\frac{1}{2}\sum_{(i,j)
\in \hat E} (x_i-x_j)^2:~\mathcal{V}(x) = \mathcal{V} \right\} =
\frac{n^2\mathcal{V}^2}{8}.$$ The last equality holds as $V(x)$,
for fixed $\mathcal{V}(x)$, is maximum when is maximum the number
of couples of elements of $x$ whose difference is equal to
$\mathcal{V}(x)$. On the other hand, denote by $\tilde E =
\{(i,i+1): ~i,i+1 \in \Gamma\}$ the edgeset of a chain network
induced by vertices in $\Gamma$. Then
$$V(x) \geq \min_{x \in \mathbb{R}^n} \left\{\frac{1}{2}\sum_{(i,j)
\in \tilde E} (x_i-x_j)^2:~\mathcal{V}(x) = \mathcal{V} \right\} =
\frac{\mathcal{V}^2}{2(n-1)}.$$ The latter equality holds as there
surely  exists a chain network defined by $\tilde E$ such that
$\sum_{(i,j) \in \tilde E} (x_i-x_j) = \mathcal{V}$. The previous
inequality holds as $\sum_{(i,j) \in E} (x_i-x_j) \geq
\mathcal{V}$ because $E$ defines a connected network on $G$, hence
there exists a path on $E$ from the agent with the maximum value
of the state and the agent with the minimum value of the state.

We can now affirm that $\frac{\mathcal{V}^2(\hat x)}{2(n-1)} \leq
V(\hat x) \leq \frac{n^2\mathcal{V}^2(\hat x)}{8}$ and that we
have $\mathcal{V}(\bar x) < {\gamma} \mathcal{V}(\hat x)$ if
$V(\bar x) \leq \frac{(\gamma\mathcal{V}(\hat x))^2}{2(n-1)}$. The
latter situation certainly occurs when $V(\bar x) \leq
\frac{4\gamma^2V(\hat x)^2}{n^2(n-1)}$.

 \qed}

\begin{theorem}\label{th:ConvergToPdelta}
Given the system (\ref{eq:xlud}) on $G = (\Gamma,E)$, implement a
distributed and stationary protocol $u(.)$ whose components have
the feedback form~(\ref{eq:lp}). For each $\nu > 0$, $0 <\gamma <
1$ there exists a finite $q(\nu,\gamma) > 0$  such that the values
assumed by the state trajectory $x(t)$ satisfy the following
condition: either $\mathcal{V}(x(t+q(\nu,\gamma)))<
{\gamma}\mathcal{V}(x(t))$ or $x(t+q(\nu,\gamma)) \in P(D,E) +
\nu$, for any $t \geq 0$ and for any disturbance realization
$\{d(t)\in D: t\geq 0\}$.

\end{theorem}

\proof{ We know that, when we implement a distributed and
stationary protocol $u(.)$ whose components have the feedback
form~(\ref{eq:lp}), the system trajectory $x(t)$ converges to a
point in $P(D,E)$ for any $x(0) \in \mathbb{R}^n$. We also know
that $V(x(t))$ is strictly decreasing along $x(t)$. Then, for any
$x(t) \in \mathbb{R}^n$, there exists a finite $\hat q(x(t),\nu,
\gamma) \geq 0$ such that either $\mathcal{V}(x(t+\hat
q(x(t),\nu,\gamma)))<\gamma\mathcal{V}(x(t))$ or $x(t+\hat
q(x(t),\nu,\gamma)) \in P(D,E) + \nu$, for any $t \geq 0$ and for
any disturbance realization $\{d(t)\in D: t\geq 0\}$. As the
system converges (exponentially) for any $x(0) \in \mathbb{R}^n$,
the value $\hat q(x,\nu,\gamma)$ is finitely bounded for $\|x\|
\rightarrow \infty$. Hence, the theorem is proved by defining
$q(\nu,\gamma) = \max_{x \in \mathbb{R}^n}\{\hat
q(x,\nu,\gamma)\}$.

\qed}

The above result is strictly related to the $\epsilon$-consensus
problem stated at the beginning of the paper (see Problem
\ref{prob:eps-consensus}). Actually, the above result means that
the state converges in finite time to a tube of radius $\epsilon=
\max\{{\cal V}(x):\,x\in P(D,E) + \nu\}$.

\section{Switching Topology.}\label{sec:Switching}
In the following, we generalize the results obtained in Section
\ref{sec:NonSwitching} to networks with switching topologies.
Consider a network $G_{\sigma(t)} = (\Gamma,E_{\sigma(t)})$ that
has a time variant edgeset $E_{\sigma(t)} \in \mathcal{E}$. We
define an edgeset $E_k$ as recurrent for a given realization of
$\sigma(t)$ if for all $t \geq 0$, there exists $t_k \geq t$ such
that $\sigma(t_k) = k$. As $\mathcal{E}$ is finite there exists at
least a recurrent edgeset for any realization of  $\sigma(t)$.
Throughout the rest of the paper we assume that all the edegsets
in $\mathcal{E}$ are recurrent over time for any realization of
$\sigma(t)$, more formally
\begin{assumption}\label{ass:RicurrentE}
Any realization $\sigma(t)$ is such that, for all $t \geq 0$,
there exists $t_k \geq t$ such that $\sigma(t_k) = k$ for all $k
\in \mathcal{I}$.
\end{assumption}
When we say that the edgeset $E_{\sigma(t)}$ is time variant we
understand that Assumption \ref{ass:RicurrentE} holds and all the
edegsets in $\mathcal{E}$ are recurrent. Note that there is no
loss of generality in Assumption \ref{ass:RicurrentE}, if it were
false, the results of this section would hold for the subset
$\hat{\mathcal{E}} \in \mathcal{E}$ of recurrent edgesets.

A basic observations is that even in presence of switches
$\mathcal{V}(x(t))$ is not increasing on $t$ as stated
in~(\ref{eq:dec}). To see this note that the protocol
$u_{\sigma(t)}(.)$ induces a continuous and bounded state
trajectory even for a network with switching topology
(\ref{eq:xlud}) on $G_{\sigma(t)} = (\Gamma,E_{\sigma(t)})$. Then,
let $t_{k}< t_{k+1}$ be two generic consecutive switching times.
Corollary~\ref{lem:initial state2} applied at time $t = t_k$
instead of $t = 0$ guarantees that $x_{i_n}(t_{k}) \leq x_i(t)
\leq x_{i_1}(t_{k})$ for all $i \in \Gamma$ and for all $t_{k} < t
< t_{k+1}$. As $x(t)$ is bounded, then $u_{\sigma(t)}(.)$ is also
bounded in the same time interval. Hence the first condition in
(\ref{eq:lp}) implies that $x(t)$ is continuous for $t_{k} < t <
t_{k+1}$, whereas, the second conditions in (\ref{eq:lp}) imposes
the continuity of the state trajectory in $t_{k}$ as $x(t_{k}^-) =
x(t_{k}^+)$ and in $t_{k+1}$ as $x(t_{k+1}^-) = x(t_{k+1}^+)$. As
a consequence, for all $t \geq 0$, we also have that $x_{i_n}(0)
\leq x_i(t) \leq x_{i_1}(0)$ for all $i \in \Gamma$ and
$\mathcal{V}(x(t))$ is not increasing.

We also need to redefine the value $\mu(.)$ initially introduced
as (\ref{eq:mu}). In particular, for given realizations $d(t)$ and
$\sigma(t)$, and a subset $Q$ of $D$, we define
\begin{equation}\label{eq:muRed} \mu(Q,E_k, t_1,t_2) = \max\{\Delta:
t_1 \leq \tilde t\leq \tilde t+\Delta \leq t_2 ~s.t.~ d(t) \in
Q~and~\sigma(t) = k~ for~all ~\tilde t\leq t \leq\tilde t+\Delta
\}.
\end{equation}
In other words, given a time interval $[t_1,t_2]$, the value
$\mu(Q,E_k,t_1,t_2)$ is the length of the longest subinterval
where $d(t)$ remains in $Q$ and $E_{\sigma(t)}$ is equal to $E_k$.

The following lemma generalizes Lemma \ref{lem:EtaUdIntersezione}
and Corollary \ref{cor:EtaUdIntersezioneMultipla}.

\begin{lemma}\label{lem:EtaUdIntersezioneRed}
Given the switched system (\ref{eq:xlud}) on $G_{\sigma(t)} =
(\Gamma,E_{\sigma(t)})$, implement a distributed and stationary
protocol $u_{\sigma(t)}(.)$ whose components have the feedback
form~(\ref{eq:lp}).  Consider a disturbance realization  $\{d(t)
\in D, ~ t\geq 0\}$ and a finite set ${\cal Q} = \{Q_1, Q_2,
\ldots\}$ of boxes of $D$. Assume that there exist two nonnegative
finite numbers $M$ and $\delta$, such that $\mu(Q_r,E_k,t,t+M) >
\delta$, for all $Q_r \in {\cal Q}$, for all $E_k \in
\mathcal{E}$, and for all $t \geq 0$. Then the equilibrium points
$x$ exist and belong to $\bigcap_{E_k \in \mathcal{E}}\bigcap_{Q_r
\in \mathcal{Q}}P(Q_r,E_k)$.
\end{lemma}
\proof{Equilibrium points exists as $\{\pi \mathbf{1}\} =
\bigcap_{E_k \in {\cal E}}\bigcap_{d \in D} P(d,E_k)$ holds. Then,
we can prove that any equilibrium point must belong to
$\bigcap_{E_k \in \mathcal{E}}\bigcap_{Q_r \in
\mathcal{Q}}P(Q_r,E_k)$ using the argument in the proof of Lemma
\ref{lem:EtaUdIntersezione} for each couple $(Q_r, E_k)$, for all
$E_k \in \mathcal{E}$ and all $Q_r \in \mathcal{Q}$.

 \qed }

The results of Lemma \ref{Lem:SignPerm} still hold in each
subinterval between two consecutive switches, and then apply even
in the switching case.

To generalize the convergence results of
Theorem~\ref{Th:EquilPoint2}, we need to introduce the following
notations, for each $Q \subseteq D$:  $L(Q,E_k) = \max_{x \in
P(Q,E_k)}\mathcal{V}(x)$ the maximum value of $\mathcal{V}(x)$ for
points in $P(Q,E_k)$; $S(Q,E_k,\nu) = \{x \in \mathbb{R}^n:
\mathcal{V}(x) \leq L(Q,E_k) +2 \nu\}$ the set of points whose
maximum difference between two components does not exceed
$L(Q,E_k) +2 \nu$. It is worth to be noted that $S(Q,E_k,\nu)$ are
tubes of radius less than or equal to $L(Q,E_k) +2 \nu$ and then,
$S(Q,E_k,\nu) \subseteq S(Q,E_{\hat k},\nu)$ whenever
$L(Q,E_k)\leq L(Q,E_{\hat k})$. Also, observe that, by definition,
it holds that $P(Q,E_k)\subseteq P(Q,E_k)+ \nu \subseteq
S(Q,E_k,\nu)$. Finally, we introduce a minimum dwell time
$\tau(\nu,\gamma) = max_{k \in \mathcal{I}}\{q_k(\nu,\gamma)\}$.
In other words, the minimum length of the switching intervals is
equal to the maximal value over the different $E_k \in
\mathcal{E}$ of the times $q(\nu,\gamma)$ introduced by
Theorem~\ref{th:ConvergToPdelta}.

\begin{theorem}\label{Th:EquilPoint4}
Given the switched system (\ref{eq:xlud}) on $G_{\sigma(t)} =
(\Gamma,E_{\sigma(t)})$, with a minimum dwell time $\hat
\tau(\nu,\gamma)$, implement a distributed and stationary protocol
$u_{\sigma(t)}(.)$ whose components have the feedback
form~(\ref{eq:lp}) with $\tilde y^{(i)}$  as in
(\ref{eq:tildeyj-1}). Assume that two values $\nu
> 0$ and $0< \gamma < 1$ are also given. The state is driven in
finite time to the tube~\begin{equation}\label{eq:T}T=\bigcap_{E_k
\in \mathcal{E}}S(D,E_k,\nu).\end{equation}
\end{theorem}

\proof{ We already know that the system trajectory $x(t)$ is
continuous and that $\mathcal{V}(x(t))$ is not increasing. Denote
by $t_s$ and $t_{s+1}$ with $t_{s+1} \geq t_s + \tau(\nu,\gamma)$
two generic consecutive switching times such that $\sigma(t)=k$
for all $t_s \leq t < t_{s+1}$. From Theorem
\ref{th:ConvergToPdelta}, we deduce that either
$\mathcal{V}(x(t_{s+1})) \leq \gamma \mathcal{V}(x(t_{s+1}))$ or
$x(t_{s+1}) \in P(D,E_k)+ \nu$. The fact that $\mathcal{V}(x(t))$
is not increasing and all the edgesets $E_k \in \mathcal{E}$ are
recurrent implies that there exists $\bar t \geq 0$ such that, if
$t_{s} \geq \bar t$, we have $x(t_{s+1}) \in P(D,E_k)+ \nu$.
Again, a not increasing $\mathcal{V}(x(t))$ implies that $x(t) \in
S(D,E_k,\nu)$, for all $t\geq t_{s+1} \geq \bar t$. As we can
apply the above argument for all the edgesets $E_k \in
\mathcal{E}$, the theorem thesis follows.

\qed}

This last result gives an answer to the $\epsilon$-consensus
problem stated at the beginning of this paper (see Problem
\ref{prob:eps-consensus}). Indeed, convergence to $T$, as in
(\ref{eq:T}), means that the agents have reached
$\epsilon$-consensus with $\epsilon = \min\{\epsilon_k\}$ where
$\epsilon_k$ is the radius of tubes $S(D,E_k,\nu)$.

Also note that the above theorem does not guarantee the
convergence to an equilibrium point. Actually, switching systems
may oscillate as shown by the following example.
\begin{example}
Consider a family of chain networks $G$ on the set of agents
$\Gamma = \{1,2,3\}$ and edgsets $E_1 = \{(1,2),(2,3)\}$ and $E_2
= \{(1,2),(1,3)\}$ (see Fig.~\ref{fig:Chain3agents}). Let $x_1(0)
= 2$, $x_2(0) = 0$, and $x_3(0) = 1$, and $\{d(t) = d=\mbox{const}
\in D: t \geq 0\}$, with in particular $d_{12}= d_{13} = d{31} =
1$ and $d_{21}= d_{23} = d_{32} = -1$. If the switching time
intervals are sufficiently long the systems trajectory oscillates
in $\mathbb{R}^3$ along the segment delimited by points $[2,0,0]$
and $[2,0,2]$. Note that only the state of agent $3$ changes over
time.
\begin{figure}[H]
    \centering%
    \includegraphics[scale=1]{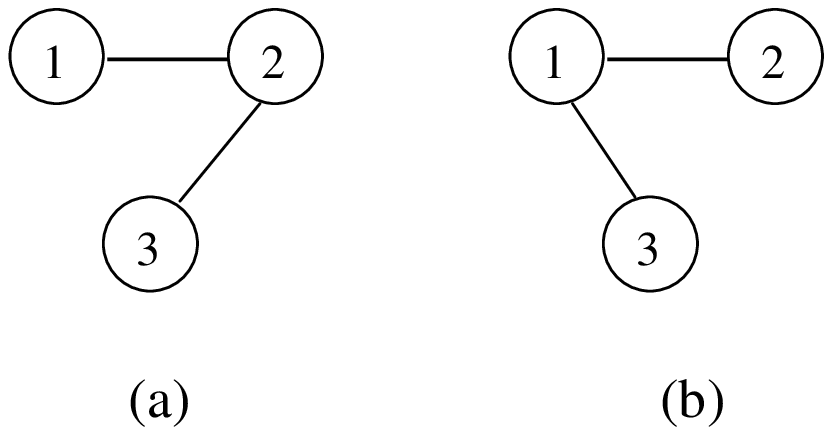}
    \caption{Family of chain networks of three agents.}
    \label{fig:Chain3agents}
\end{figure}

\end{example}

Finally, note that we could generalize
Theorem~\ref{Th:EquilPoint3} and
Corollary~\ref{cor:TheoEtaUdIntersezioneMultipla} only in the
assumption that, for all $Q_r \in {\cal Q}$, for all $E_k \in
\mathcal{E}$, the values $\mu(Q_r,E_k,t,t+M)$ define sufficiently
long intervals $[t_s, t_{s+1})$ so that either
$\mathcal{V}(x(t_{s+1}))$ is finitely smaller than
$\mathcal{V}(x(t_{s+1}))$ of a finite value or $x(t_{s+1}) \in
P(Q_r,E_k) + \nu$. In this case we have the system trajectory
eventually assume values in~$\bigcap_{E_k \in
\mathcal{E}}\bigcap_{Q_r \in \mathcal{Q}}S(Q_r,E_k,\nu)$.

%
%

\section{Conclusions.}\label{sec:conclusions}
Despite the literature on consensus is now becoming extensive,
only  few approaches have considered a disturbance affecting the
measurements. In our approach we have assumed an UBB noise in the
neighbors' state feedback as it requires the least amount of
a-priori knowledge on the disturbance. Only the knowledge of a
bound on the realization is assumed, and no statistical properties
need  to be satisfied. Because of the presence of UBB disturbances
convergence to equilibria with all equal components is, in
general, not possible. Therefore, the main contribution has been
the introduction and solution of the $\epsilon$-consensus problem,
where the states converge in a tube of ray $\epsilon$
asymptotically or in finite time. In solving the
$\epsilon$-consensus problem we have focused on linear protocols
and presented a rule for estimating the average from a compact set
of candidate points.

\end{document}